\documentclass[12pt,french,brochure]{bettyb}

\usepackage{url}

\usepackage{natbib}
  
\makeatletter
  \def\NAT@nmfmt#1{{\scshape\NAT@up#1}}
  \def\bibliography#1{%
    \if@filesw
      \immediate\write\@auxout{\string\bibdata{#1}}%
    \fi
    \expandafter\input{\bbl@main@language bst.tex}%
    \@input@{\jobname.bbl}        
  }
\makeatother

\iffalse
  \usepackage{fontspec}
  \setromanfont[Ligatures=TeX]{TeX Gyre Pagella}
  \usepackage{unicode-math}
\else
  \usepackage[utf8]{inputenc}
  \usepackage{newpxtext,newpxmath}
\fi

\usepackage[T1]{fontenc}
\usepackage[french]{babel}
\let\mathscr \mathcal

\usepackage{hyperref}

\author{Antoine Chambert-Loir}
\address{Université Paris Cité \\ 
IMJ-PRG \\
F-75013, Paris, France}
\email{antoine.chambert-loir@u-paris.fr}
\title{Les conjectures de Weil : \\ origines, approches, généralisations}

\date{Novembre 2022}
\bbkannee{6\textsuperscript{e} ann\'ee, 2022--2023}
\bbknumero{33}

\begin{abstract}
Je retracerai l’histoire des conjectures de Weil 
sur le nombre de solutions d’équations
polynomiales dans un corps fini
et quelques unes des approches qui ont été proposées
pour les résoudre.
\end{abstract}

\alttitle{The Weil conjectures: origins, approaches, generalizations}
\begin{altabstract}
I recount the history of the conjectures by Weil on the number
of solutions of polynomial equations in finite fields,
and some of the approaches that have been proposed to solve them.
\end{altabstract}

\def\Q{\mathbf Q}
\def\A{\mathbf A}
\def\R{\mathbf R}
\def\C{\mathbf C}
\def\Z{\mathbf Z}
\def\P{\mathbf P}
\def\F{\mathbf F}
\let\eps\varepsilon

\def\div{\operatorname{div}}
\def\Tr{\operatorname{Tr}}
\def\Div{\operatorname{Div}}

\def\SL{\operatorname{SL}}

\def\Card{\operatorname{Card}}

\def\abs#1{\left\lvert#1\right\rvert}

\newcommand{\legendre}[2]{\genfrac{(}{)}{1pt}{}{#1}{#2}}

\begin{document}
\maketitle

\section{Prologue: Gauss}

Dans ses \emph{Disquisitiones arithmeticae}, \cite{Gauss-1863a}
démontrait plusieurs théorèmes qui dénombrent les solutions
de certaines équations en congruences modulo un nombre premier.
Il prouve par exemple que le nombre de couples $(x,y)$
d'entiers modulo~$p$ tels que $x^2+y^2\equiv 1 \pmod p$
est donné par 
\[ N (x^2+y^2=1) = p - \legendre{-1}p, \]
où $\legendre \cdot p$ désigne le « symbole de Legendre » modulo~$p$,
défini par $\legendre ap=0$ si $p$ divise~$a$,
$1$ si $a$ est un carré modulo~$p$ (l'expression était « résidu quadratique »),
et $-1$ sinon.
Il utilisait ensuite cette formule pour établir la « loi complémentaire »
de sa \emph{loi de réciprocité quadratique} (theorema aureum)
disant que $\legendre 2p$ vaut $1$ si $p\equiv \pm1\pmod 8$
et $-1$ si $p\equiv \pm3\pmod 8$.

La dernière entrée de son agenda, publié par \cite{Klein-1903},
est une affirmation du même genre:
\begin{quote}
Observatio per inductionem facta gravissima theoriam residuorum
biquadraticorum cum functionibus lemniscaticis elegantissime nectens.
Puta si $a+bi$ est numerus primus, $a-1+bi$ per $2+2i$ divisibilis, 
multitudo omnium solutionum congruentiae
\[ 1 \equiv xx+yy+xxyy \pmod{a+bi} \]
inclusis
\[ x = \infty, \quad y = \pm i; \qquad x = \pm i, \quad y=\infty\]
fit 
\[ = (a- 1)^2+ bb. \]
\leftskip0pt plus 1fil \parfillskip0pt  1814 Iul. 9. 
\end{quote}

Autrement dit, on considère un nombre premier~$p$ congru
à 1 modulo~4. D'après Fermat, il s'écrit
sous la forme $a^2+b^2$ (de sorte que $\pi=a+bi$ est un premier
de l'anneau $\Z[i]$ des entiers de Gauss).
Quitte à échanger~$a$ et~$b$, on suppose que $a$ est impair et $b$ est pair ;
quitte à remplacer~$a$ par~$-a$,  on peut supposer que $a+b\equiv 1\pmod 4$;
alors, $a-1+bi$ est divisible par~$2+2i$ dans~$\Z[i]$.
De manière équivalente, on impose
que $(a-1)^2+b^2$ est divisible par~$8$,
d'où $p-2a+1\equiv 0 \pmod 8$;
lorsque $p\equiv 1\pmod 8$, on a donc $a\equiv 1\pmod 4$,
tandis que lorsque $p\equiv 5\pmod 8$, on a $a\equiv 3\pmod 4$.
Gauss affirme alors que le nombre de couples $(x,y)$
dans~$\Z$ tels que $x^2+y^2+x^2y^2=1$ modulo~$p$
est égal à $(a-1)^2+b^2-4=p-2a-3$.

C'est toutefois une \emph{conjecture}
que Gauss énonce ici ---
\emph{observatio per inductionem facta gravissima} ---
et plusieurs mathématiciens après lui en proposeront
des démonstrations: \cite{Herglotz-1921,Chowla-1949}.

Pourtant, ainsi que le rappelle~\cite{Weil-1949}, 
Gauss avait démontré des résultats analogues
pour les congruences cubiques
\[ ax^3-by^3\equiv 1 \pmod p \]
(\citet[\S358]{Gauss-1863b})
et biquadratiques
\[ a x^4-by^4\equiv 1 \pmod p,
\qquad y^2\equiv ax^4-b \pmod p \]
(\citet[\S23]{Gauss-1863a})
On a par exemple
\[ N(y^2=1-x^4) = \begin{cases}
2 & \text{si $p=2$;} \\
 p - 1 & \text{si $p\equiv 3 \pmod 4$;}  \\
 p - 1 - 2a & \text{si $p\equiv 1\pmod 4$, } \end{cases} \]
où l'on a écrit $p=a^2+b^2$ avec $a\equiv p\pmod 8$.

\section{Estimation et dénombrement}

\cite{HasseDavenport-1935} et \cite{Weil-1949} généralisent 
ces formules à toutes les équations de la forme
\[ \sum_{i=1}^m a_i x_i^{n_i} = 1 \]
non seulement en congruences modulo un nombre premier~$p$, 
c'est-à-dire dans le corps fini~$\Z/p\Z$,
mais dans un corps fini arbitraire.
Leur solution est exprimée en termes de sommes de Gauss 
ou de Jacobi (introduites par Gauss) et met en évidence
de très grandes régularités.

Ces questions, et notamment la réflexion de Weil,
sont nées de la fusion de deux motivations assez différentes.

Dans sa thèse, \cite{Artin-1924a,Artin-1924} avait poursuivi
l'analogie entre corps de nombres et corps de fonctions 
entamée par \cite{Dedekind-1857} en mettant en vis-à-vis
le corps~$\Q$ des nombres rationnels et les corps de
fractions rationnelles $\F_p(x)$ 
en une indéterminée~$x$ à coefficients dans~$\Z/p\Z$.
Il met en regard l'anneau~$\Z$ des entiers relatifs et celui~$\F_p[x]$
des polynômes, tous deux des anneaux principaux, 
les nombres premiers correspondent aux polynômes irréductibles unitaires,
les unités~$\pm1$ aux polynômes constants non nuls et
les entiers strictement positifs aux polynômes unitaires.
En observant que pour $n>0$, on a $n=\Card(\Z/(n))$
tandis que pour un polynôme unitaire~$P$, on a $\Card(\F_p[x]/(P))=p^{\deg(P)}$,
l'analogue de la fonction zêta de Riemann
\[ \zeta_{\Z}(s) = \sum_{0 < n } \frac1{n^s} = \prod_{\text{$p$ premier}}
 \frac1{1-p^{-s}} \]
est la série 
\[ \zeta_{\F_p[x]} (s) = \sum_{0\neq I \subseteq \F_p[x]} \frac1{\Card(\F_p[x]/I)^{s}}
 = \prod_{\text{$P$ irréductible unitaire}} \frac1{1- p^{s \deg(P)}} .\]

Artin poursuivait donc cette analogie 
des nombres rationnels aux corps quadratiques 
 en étudiant l'arithmétique
des corps de la forme 
$ \F_p(x) (\sqrt{f(x)}) $ obtenus
en adjoignant à~$\F_p(x)$ une racine carrée d'un
polynôme sans facteur carré~$f$.
Il introduisait en particulier leur fonction zêta 
et démontrait que c'est une fraction rationnelle de $q^{-s}$,
établissait leur équation fonctionnelle, 
vérifiait la « formule analytique du nombre de classes »,
les utilisait pour étudier l'analogue des théorèmes des nombres
premiers et de la progression arithmétique.
Il vérifiait aussi, mais seulement dans un petite nombre
de cas (\S23),
l'analogue de l'\emph{hypothèse de Riemann},
c'est-à-dire que les zéros de ces fonctions zêta
soient des nombres complexes de partie réelle~$1/2$;
il avait apparemment en vue la finitude de l'ensemble
des corps de ce type dont le nombre de classes est donné.
En des termes géométriques encore anachroniques,
dont la nécessité n'apparaîtra que peu à peu,
il étudie la courbe affine hyperelliptique d'équation $y^2=f(x)$.

En Grande-Bretagne, \cite{Davenport-1933},
alors élève de L.~J. Mordell,
tentait d'estimer des « sommes exponentielles »
analogues aux sommes de Gauss, 
du genre
\[ \sum_{x\in \Z/p\Z} \legendre {f(x)}p \]
où $f$ est un polynôme.  C'est une somme de $p$ termes
égaux à~$\pm1$, et parfois~$0$; guidé par l'heuristique 
qu'elle devrait être de l'ordre de~$\sqrt p$, 
il obtenait une majoration non triviale, en $\mathrm O(p^{3/4})$
lorsque $f$ est unitaire de degré~$4$.
En observant que $ 1 + \legendre{f(x)}p $ est le nombre
de solutions dans~$\Z/p\Z$ de l'équation $y^2=f(x)$,
on a 
\[ \sum_{x\in \Z/p\Z} \legendre {f(x)}p = N(y^2=f(x)) - p \]
faisant un lien entre les deux questions.

Peu après sa thèse, Artin avait étendu son étude
en remplaçant $\Z/p\Z$ par un corps fini arbitraire.
Ainsi, il a pu montrer que la validité de son hypothèse de Riemann
était inchangée si l'on étendait le corps fini, sans modifier
le polynôme~$f$ et en déduire des familles de polynômes
pour lesquelles l'hypothèse est vérifiée. 
Selon \cite{Roquette-2018}, il semble que ce soit
cette vérification qui lui ait permis de croire en la véracité
de l'hypothèse de Riemann. Roquette publie cependant
une lettre qu'Artin écrit à Herglotz
dans laquelle il se  plaint de l'attitude 
quelque peu « mandarinale » de Hilbert
lors de son exposé à Göttingen et explique 
qu'il souhaite subitement abandonner ce sujet. 
(Pour une description de ces lettres 
d'Artin à Herglotz et les travaux non publiés
d'Artin, voir aussi \cite{Ullrich-2000a}.)
Dans les années suivantes,
M.~Deuring et F.K.~Schmidt ont mis en place la théorie des courbes
sur un corps fini, notamment leurs fonctions zêta.
C'est cependant Artin qui a expliqué à H.~Hasse,
lors d'une visite à Hamburg en 1932, le lien entre le problème
de majoration de sommes d'exponentielles et 
l'analogue de l'hypothèse de Riemann, convaincant le second de s'y attaquer:
en 1935, Hasse résoudra  le cas des courbes elliptiques,
c'est-à-dire lorsque le polynôme~$f$ a degré~$3$ ou~$4$.

\section{Les conjectures de Weil}

Pour énoncer les conjectures générales de Weil,  
nous adoptons désormais un langage géométrique.
Considérons donc un corps fini~$k$
et une variété algébrique~$V$ sur~$k$, c'est-à-dire, suivant
son goût, un schéma de type fini sur~$k$,
ou bien un système d'équations polynomiales 
$\{ f_1(x_1,\dots,x_m)=\dots = f_r(x_1,\dots,x_m)=0\}$ 
à coefficients dans~$k$.

L'ensemble $V(k)$ des points $k$-rationnels de~$V$ 
correspond aux solutions dans~$k^m$ du système $\{f_1=\dots=f_m=0\}$;
comme le corps~$k$ est fini, c'est un ensemble fini
et on note $N(V)$ son cardinal.

Soit $q$ le cardinal de~$k$.
La théorie des corps finis apprend que $q$
est une puissance de la caractéristique~$p$ de~$k$
et qu'inversement, toute puissance de~$p$ est le cardinal
d'un corps fini, unique à isomorphisme près. Ainsi,
toute puissance~$q^n$ de~$q$ est le cardinal d'un corps fini~$k_n$
qui est l'unique extension de degré~$n$ de~$k$.
On notera $N_n(V)$ le cardinal de $V(k_n)$.

Introduisant une indéterminée~$t$, \citet{Weil-1949} 
définit alors la série formelle à coefficients rationnels
\[ Z_V(t)  = \exp \big( \sum_{n=1}^\infty \tfrac1n N_n(V) t^n \big). \]
Motivée par les calculs que Weil avait faits au début de son article, 
cette formule peut-être étrange généralise les définitions de Dedekind et Artin.
En effet, les points fermés~$x$ du schéma~$V$ correspondent
aux idéaux maximaux~$M$ de l'anneau $A=k[x_1,\dots,x_m]/(f_1,\dots,f_r)$,
ou aux solutions $(x_1,\dots,x_m)$ 
du système $\{f_1,\dots,f_r\}$ dans une clôture algébrique
fixée~$\overline k$ de~$k$;
le corps résiduel $\kappa(x)$ correspond au corps résiduel~$A/M$
et aussi à l'extension finie de~$k$ engendrée par~$x_1,\dots,x_m$;
on note $\deg(x)$ le degré de cette extension.
Alors, on a
\[ Z_V(t) = \prod_{x \in \abs V} \frac1{1-t^{\deg(x)}},  \]
une formule qui prouve que $Z_V(q^{-s})=\zeta_V (s)$
dans les cas considérés par Dedekind et Artin.
Elle montre aussi que $Z_V(t)$ est à coefficients entiers

Dans le cas de Dedekind, $V=\A^1$ est la droite affine,
on a $N_n(V)=q^n$ et on obtient 
\[ Z_{\A^1}(t) = \exp\big( \sum_{n=1}^\infty \tfrac1n q^n t^n\big)
 = \dfrac{1}{1-qt}. \]
Les cas de l'espace affine et de l'espace projectif sont également
élémentaires, on trouve
\[ Z_{\A^d}(t) = \frac1{1-q^d t} \quad\text{et}\quad 
 Z_{\P_d}(t) = \frac1{(1-t)(1-qt)\dots (1-q^d t)}. \]
\citet{Weil-1949} considère aussi le cas des grassmanniennes; il observe 
dans tous ces cas, ainsi que celui
des hypersurfaces « diagonales »
d'équation $ a_1 x_1^{n_1}+\dots+a_m x_m^{n_m}=1$
(ou leur version homogène),
un lien entre le résultat obtenu
et le polynôme de Poincaré de l'espace affine, l'espace projectif,
la grassmannienne ou l'hypersurface diagonale complexe analogue.

La première conjecture de Weil énonce que
\emph{$Z_V(t)$ est une fraction rationnelle.}

L'analogie avec la situation géométrique classique peut même être
poursuivie. Pour cela, il convient de faire l'hypothèse supplémentaire
que $V$ est un schéma propre, lisse et géométriquement connexe;
notons~$d$ sa dimension.
Sa compagne en géométrie complexe est alors une variété différentielle
compacte et connexe de dimension (réelle)~$2d$.
Weil suggère en effet qu'on peut écrire $Z_V$ sous la forme
\[ Z_V(t) = \frac{P_1(t) \dots P_{2d-1}(t)}
{P_0(t)\dots P_{2d}(t)} \]
où $P_0,\dots,P_{2d}$ sont des polynômes 
de terme constant~$1$ et de degré les \emph{nombres de Betti}
$b_0,\dots,b_{2d}$ de la compagne complexe de~$V$.

Les seconde et troisième conjectures se placent encore sous
cette hypothèse.

Weil postule l'existence d'une \emph{équation fonctionnelle} du type
\[ Z_V(1/q^d t) = \eps q^{\chi d/2}t^{\chi} Z_V(t), \]
où $\eps=\pm1$ et $\chi=b_0-b_1+b_2-\dots$ 
est la caractéristique d'Euler-Poincaré de~$V$.

Enfin, la troisième conjecture affirme que les
inverses des racines de~$P_i$ sont des entiers algébriques
de valeur absolue $q^{i/2}$.

Dans son exposé \citep{Weil-1956a}
au Congrès international d'Amsterdam,
il précisera son intuition géométrique, dans le cas des courbes,
d'une façon qui orientera définitivement les études
ultérieures.

\section{Cohomologies de Weil}

Considérons un schéma propre, lisse, géométriquement intègre~$V$ 
sur un corps fini~$k$ de cardinal~$q$.
Dans le contexte des courbes elliptiques,
Hasse avait introduit le morphisme de Frobenius de~$V$.
Ici, c'est un morphisme de schémas $\phi_V\colon V\to V$ 
qui est l'identité sur les points
mais, du point de vue des anneaux locaux de~$V$,
est donné par l'élévation à la puissance~$q$.
Chaque ouvert affine~$U$ de~$V$ est stable par~$\phi_V$,
et sur l'anneau~$A$ de~$U$, $\phi_V$ est donné 
par le morphisme de $k$-algèbres $a\mapsto a^q$.

Si l'on étend les scalaires de~$k$ à une clôture algébrique~$\overline k$
de~$k$, le morphisme~$\phi_V$ devient un morphisme de schémas
$\phi_{\overline V} \colon \overline V\to {\overline V}$
et $V(k)$ apparaît comme l'ensemble des \emph{points fixes}
de $\phi_{\overline V}$ agissant sur~$\overline V$.

Comme l'explique~\citet{Weil-1956a},
la formule des traces de Lefschetz en topologie algébrique
suggère alors une expression
\[ N(V) = \sum_{i=0}^{2d} (-1)^i \Tr \big(\phi_{\overline V}^*
\mid H_i(\overline V)\big) \]
où les $H^i(\overline V)$ seraient les groupes d'homologie
de $\overline V$, définis fonctoriellement en~$\overline V$.

Notons que dans la formule précédente, les points fixes
sont comptés avec multiplicité~1. En effet, 
comme la dérivée du polynôme $T^q$, égale à $qT^{q-1}$, est nulle 
dans~$k[T]$, la « différentielle » de $\phi_{\overline V}$
est nulle, et en particulier, son graphe dans $\overline V\times \overline V$
est transverse à la diagonale.

Appliquant cette idée aux puissances $\phi_{\overline V}^n$ de~$\phi_{\overline V}$, on obtient une formule
\[ N_n(V) = \sum_{i=0}^{2d} (-1)^i \Tr \big((\phi_{\overline V}^n)^* \mid
H_i(\overline V)\big), \]
qui conduit à l'expression
\[ Z_V(t) = \frac{P_1(t)P_3(t)\dots P_{2d-1}(t)}
{P_0(t)P_2(t)\dots P_{2d}(t)} \]
où 
\[ P_i (t)  = \det\big(1-t \phi_{\overline V}\mid H_i( \overline V)\big). \]

Dans cette analogie, l'équation fonctionnelle traduirait
la dualité de Poincaré 
$H_i(\overline V) \leftrightarrow H_{2d-i}(\overline V)$.

Ces idées suggèrent également des formes d'uniformité 
pour des familles de variétés. Un théorème d'Ehresmann
énonce qu'une submersion propre de variétés différentielles 
est localement triviale, une fibre est difféomorphe
aux fibres voisines. En géométrie algébrique, cela suggère que 
 si $\mathscr V\to S$ est
un morphisme propre et lisse de schémas, à fibres géométriquement connexes,
les espaces d'homologie $H_i(\overline{\mathscr V_s})$ devraient être
de dimension constante, lorsque $s$ parcourt~$S$.
Lorsque $\kappa(s)$ est un sous-corps de~$\C$,
la variété $\mathscr V_s$ peut être considérée comme une
variété complexe, et 
on voudrait également que les espaces
d'homologie soient  reliés à ceux donnés par la topologie.

Les années 1960 ont vu le début d'une quête 
des « cohomologies\footnote{Je ne sais pas bien comment
l'on est passé de l'homologie au point de vue dual de la cohomologie.}
de Weil ».
Les constructions de \cite{Weil-1948} éclairent le cas des courbes.
Dans le cas topologique, l'espace
d'homologie $H_1(V)$ d'une courbe~$V$ (projective, lisse,
connexe) de genre~$g$ est un $\Q$-espace vectoriel
de dimension~$2g$ et
les travaux de Weil suggèrent, lorsque $V$
est une telle courbe sur un corps~$k$ qu'on puisse
prendre pour $H_1(V)$ un espace vectoriel
de dimension~$2g$, non pas sur~$\Q$,
mais sur le corps~$\Q_\ell$ des nombres $\ell$-adiques,
où $\ell$ est un nombre premier distinct de la caractéristique de~$k$.
Chez Weil, cet espace vectoriel est
construit à partir du « module de Tate » $\ell$-adique 
de la jacobienne~$J$ de~$V$, 
c'est-à-dire la limite
\[ \varprojlim_n J[\ell^n] \]
des groupes des points d'ordre une puissance de~$\ell$ de la
variété abélienne~$J$.
Comme $V[\ell^n]$ est isomorphe $(\Z/\ell^n\Z)^{2g}$,
ce module est un $\Z_\ell$-module libre de rang~$2g$,
dont on peut prendre le produit tensoriel par~$\Q$.

Dès le début, cette quête 
est néanmoins placée sous l'ombre d'une mise en garde de Serre
quant au type d'espaces vectoriels que l'on peut espérer obtenir.  
Supposons en effet que $V$ soit une courbe elliptique \emph{supersingulière}.
L'anneau des endomorphismes de~$V$ qu'avait introduit Hasse,
est alors un ordre d'une algèbre de quaternions~$A$ sur~$\Q$
et la fonctorialité de l'homologie fournit 
un plongement de cette algèbre $A$ dans l'anneau des matrices $2\times 2$
à coefficients dans~$F$. Cela impose que $A\otimes F$
soit isomorphe à l'algèbre~$M_2(F)$, ce qui est effectivement
le cas lorsque $F=\Q_\ell$ pour $\ell$ différent de la caractéristique,  
mais n'a pas lieu pour $F=\R$ ou $F=\Q_p$.
En particulier, la quête d'une cohomologie de Weil
à valeurs « $\Q$-espaces vectoriels » est illusoire.

Parmi les quelques exemples de cohomologies de Weil maintenant
à notre disposition,
citons:
\begin{enumerate}
\item 
La \emph{cohomologie étale} définie par Grothendieck et développée
avec plusieurs collaborateurs, dont M.~Artin, J.-L.~Verdier
et P.~Deligne
au sein des \emph{Séminaires de géométrie algébrique}
(volumes~4, 5 et 7 en particulier). 
Sa définition repose sur la notion de revêtement en géométrie algébrique;
qu'elle donne le résultat attendu
utilise le calcul des revêtements des courbes en termes de jacobiennes,
ainsi qu'un théorème de M.~Artin selon lequel 
les variétés algébriques complexes ont une base d'ouverts
qui sont des $K(\pi,1)$, c'est-à-dire
dont le type d'homotopie est gouverné par leur groupe fondamental.

\item
Les \emph{cohomologies cristalline} et \emph{rigide} définies par Berthelot
sur le modèle d'une construction par Grothendieck de la cohomologie
de De Rham. Le modèle est la théorie des équations différentielles.
\end{enumerate}
En quelque sorte ces deux cohomologies
sont les deux pôles de la correspondance de Riemann--Hilbert.

De nombreuses questions restent cependant ouvertes:
on ne sait par exemple pas si la dimension des groupes
de cohomologie étale dépend du choix du nombre premier~$\ell$!

\section{Approches des conjectures de Weil (rationalité et équation fonctionnelle)}

Il y a essentiellement trois preuves
de la rationalité des fonctions zêta de Weil.

1) La première, limitée au cas des courbes,
est due à \cite{Schmidt-1931} et repose sur le théorème
de Riemann--Roch. Donnons-en le principe.
Lorsqu'on développe la formule exprimant $Z_V(t)$
comme un produit sur les points fermés de~$V$,
on obtient une formule 
\[ Z_V(t) = \sum_n \Card(\Div_n^+(V)) t^n, \]
où $\Div_n^+(V)$ désigne l'ensemble des diviseurs effectifs
de degré~$n$ sur~$V$, c'est-à-dire des combinaisons linéaires
formelles de points fermés dont la somme des degrés vaut~$n$.
Admettons qu'il existe un diviseur de degré~1, $D_0$, sur~$V$
(\citet{Schmidt-1931} le déduit de l'analyse ci-dessous;
sinon, la fonction zêta de~$V$ sur une extension convenable
de~$k$ n'aurait pas un pôle simple en $1$ et $1/q$)
et associons à un diviseur~$D\in\Div_n^+(V)$
la classe $[D-nD_0]$ de $D-nD_0$ dans la jacobienne de~$V$.
Pour tout diviseur~$E$ de degré~$0$, on trouve
les diviseurs effectifs $D$ tels que $[D-nD_0]=E$
en considérant l'espace de Riemann--Roch $\mathscr L(E+nD_0)$:
les fonctions rationnelles sur~$V$ dont le diviseur~$f$ 
vérifie $\div(f)+E+nD_0\geq 0$; alors $D=\div(f)+E+nD_0$ 
est un diviseur effectif de degré~$n$ sur~$V$ et $D-nD_0$ a
même classe que~$E$.
Cela permet de calculer le cardinal des éléments 
de $\Div_n^+(V)$ de classe~$[E]$:
\[ \frac{q^{h(E+nD_0)}-1}{q-1}, \]
où $h(E+nD_0)=\dim (\mathscr L(E+nD_0))$.
D'après le théorème de Riemann--Roch, cette dimension vaut
$n+1-g$ si $n>2g-2$; en ajoutant les contributions
de chaque classe de diviseur~$E$ de degré~$0$, on obtient
que $(1-t)(1-qt)Z_V(t)$ est un polynôme de degré~$2g$.

Le théorème de Riemann--Roch affirme plus précisément
que $h(E+nD_0)=n+1-g +h(K_V-E-nD_0)$, où $K_V$ est un diviseur
canonique sur~$V$. En termes plus modernes, 
le théorème de Riemann--Roch calcule
une caractéristique d'Euler-Poincaré
$h^0(E+nD_0)-h^1(E+nD_0)=n+1-g$,
on a $h=h^0$ et la « dualité de Serre »
exprime $h^1(E+nD_0)=h^0(K_V-E-nD_0)$.

2) \citet{Dwork-1960a} a démontré la rationalité des fonctions
zêta en toute dimension par une méthode très originale.
Il se ramène d'abord au cas où $V$ est une hypersurface
d'équation $f=0$ dans un tore $(\A^1\setminus\{0\}) ^m $
et récrit $Z_V(t)$ comme  une somme exponentielle.
Par des arguments d'analyse $p$-adique, il prouve que cette
série est le quotient de
deux séries entières dont le rayon de convergence $p$-adique
est infini; en particulier, $Z_V(t)$ apparaît comme
une \emph{fonction méromorphe} définie sur tout $\Z_p$.
En utilisant que $\Card(V(k_n))\leq q^{nm}$, on constate
que cette même série définit une fonction holomorphe
sur le disque ouvert de rayon $q^{-m}$ dans~$\C$.
En adaptant le critère de rationalité de \cite{Borel-1894},
Dwork en déduit que $Z_V(t)$ est le développement
de Taylor d'une fraction rationnelle.

3) Comme esquissé au paragraphe
précédent, les cohomologies de Weil fournissent une démonstration
de la rationalité des fonctions zêta. Il y a en particulier
une démonstration par voie étale (SGA~5) et une démonstration
par voie cristalline.

L'équation fonctionnelle des fonctions zêta
se déduit alors de la dualité de Poincaré.

\section{Approches des conjectures de Weil (hypothèse de Riemann)}

\subsection{Le cas des courbes}

Au delà de quelques familles d'exemples dans l'article
initial d'\cite{Artin-1924} et dans l'étude
qu'il avait abandonnée à l'automne~1921,
la première démonstration concerne le cas des courbes elliptiques
pour lesquelles Hasse a offert  deux preuves.

La première reposait sur la théorie de la multiplication complexe,
à savoir la construction d'une courbe elliptique complexe
dont les endomorphismes soient un ordre d'un corps quadratique
imaginaire et telle qu'en en « réduisant modulo~$p$ » l'équation
de Weierstrass, on obtienne la courbe elliptique initiale~$V$ sur~$k$.
L'endomorphisme de Frobenius de~$V$ correspond alors
à un nombre quadratique imaginaire~$\pi$ tel que $\pi\overline\pi=q$
et le cardinal de~$V(k)$ est égal à $\abs{1-\pi}^2=1-\Tr(\pi)+q$.
D'une certaine manière, cette preuve est dans la lignée des
calculs de Gauss: la courbe d'équation $y^2=1-x^4$
est en effet une courbe elliptique à multiplication complexe par~$\Z[i]$.

La seconde démonstration \citep{Hasse-1936} repose sur une étude
plus algébrique de l'anneau des endomorphismes de~$V$.
Hasse démontre que le degré d'un tel endomorphisme
fournit une forme quadratique définie positive
et les égalités $N(V)=\deg(1-\phi_V)$ et $\deg(\phi_V)=q$,
puis en déduit l'inégalité
$\abs{N(V)-1-q}\leq 2\sqrt q$.

La démonstration générale est due à Weil.
L'anneau des endomorphismes d'une courbe elliptique
est remplacé par celui des \emph{correspondances}
qu'avait introduit Deuring. Weil démontre
une formule $N(V)=1-\Tr(\phi_V)+q$, où $\Tr(\phi_V)$ 
désigne la trace de l'action de~$\phi_V$ sur ce qu'on appelle
aujourd'hui le module de Tate $\ell$-adique de la jacobienne de~$V$.
De fait, la note \citep{Weil-1940} établissait cette formule
sous l'hypothèse que la jacobienne de~$V$ possède, pour tout entier~$n$
premier à la caractéristique de~$k$, précisément $n^{2g}$
points~$\alpha$ tels que $n\alpha=0$,
de sorte que ce module de Tate est de rang~$2g$ sur l'anneau 
des entiers $\ell$-adiques.
L'année suivante, \cite{Weil-1941} fournirait les détails
d'une preuve de cette formule.

Presque dix années furent nécessaire à la rédaction
détaillée de cette preuve, pour laquelle Weil dut  
développe une « géométrie algébrique abstraite »,
sur un corps arbitraire et en particulier 
indépendante de toute considération topologique
\citep{Weil-1946,Weil-1948,Weil-1948b}.

Le cœur de l'hypothèse de Riemann apparaît comme la positivité
de la trace de certaines correspondances
sur la courbe~$V$. Par définition, une telle correspondance
est une somme formelle de courbes tracées 
sur la surface~$V\times V$, que l'on imagine comme
le graphe d'une « fonction multivaluée » de~$V$ dans~$V$;
les correspondances se composent naturellement,
et une correpondance~$C$ possède une correspondance symétrique~$C'$,
obtenue en échangeant l'ordre des facteurs.
Elles agissent sur les diviseurs sur~$V$,
en particulier sur le module de Tate.
L'hypothèse de Riemann
apparaît comme conséquence de l'inégalité $\Tr(C\circ C')>0$,
que Weil démontre via l'étude des endomorphismes de la jacobienne de~$V$.

\cite{MattuckTate-1958} et \cite{Grothendieck-1958b}
reprouvent cette inégalité de façon plus directe, sans recours
à la jacobienne, via le théorème de Riemann--Roch pour la
surface $V\times V$. C'est le théorème de l'indice de Hodge
qui apparaît alors comme source de l'inégalité cruciale.

Dans le cas des courbes hyperelliptiques,
\cite{Stepanov-1969} a proposé une démonstration 
encore plus élémentaire, en ce qu'elle se place sur la courbe
elle-même et n'utilise que le théorème de Riemann--Roch.
Elle a été généralisée indépendamment par \cite{Schmidt-1976}
et \cite{Bombieri-1974}.
L'idée est de construire une fonction rationnelle~$f$ 
non nulle sur la courbe~$V$ de degré contrôlé 
et qui, par construction,
s'annule automatiquement en chaque point de~$V(k)$; 
il en résulte une inégalité de la forme $N(V)\leq q + 1 + \mathrm O( \sqrt q)$.
En utilisant cette inégalité pour des courbes auxiliaires
convenables~$V'$, on obtient une inégalité dans l'autre sens
$N(V)\geq q+1 - \mathrm O(\sqrt q)$. Un argument simple fondé
sur la rationalité de~$Z_V$ permet alors de conclure.

\subsection{Estimées partielles en dimension supérieure}

Pendant près de 25 ans, la seule estimée générale 
était due à \citet{LangWeil-1954}. Pour une sous-variété
$V\subset\P_n$, 
de dimension~$d$ et géométriquement irréductible,
ils démontrent une majoration du type
\[ \abs{N(V) - q^{d} } \leq c_1 q^{d-\frac12}  + c_2 q^{d-1}, \]
où $c_1=(\deg(V)-1)(\deg(V)-2)$ 
et $c_2$ est un nombre réel non explicite,
mais qui ne dépend que de~$n$ et de $\deg(V)$.

Modulo l'utilisation de l'hypothèse de Riemann pour les courbes,
leur démonstration est très simple: elle consiste à considérer
les intersections $V\cap H$ de~$V$ avec un sous-espace projectif variable~$H$
de dimension $n+1-d$. Supposons pour simplifier que $V$ ne soit pas
contenue dans un sous-espace projectif non trivial.
Par un théorème de type Lefschetz, ces intersections $V\cap H$
sont le plus souvent des courbes géométriquement irréductibles,
de genre arithmétique $g=(\deg(V)-1)(\deg(V)-2)$.
Quitte à prendre garde à leurs singularités, elles donnent toutes
lieu à une estimation du type
\[ \abs{N(V\cap H) - q }  \leq g \sqrt q + c, \]
où $c$ est un nombre réel qui ne dépend pas de~$H$.
Par des arguments de géométrie algébrique
les sous-espaces projectifs~$H$ pour lesquels $V\cap H$
n'est pas géométriquement irréductible figurent parmi
une sous-variété de la variété grassmannienne adéquate;
dans leur cas, on se contente d'une majoration $N(V\cap H)\leq \deg(V)$.
Il reste à sommer ces estimations lorsque $H$ varie;
celles du premier type dominent l'asymptotique et fournissent
le premier terme $c_1q^{d-\frac12}$, les secondes donnent
lieu au terme d'erreur.

\subsection{Conjectures standard}

Dans son exposé au Congrès international, \cite{Weil-1956a}
avait expliqué une démonstration de la positivité
$\Tr(C\circ C')>0$ qui utilise la géométrie algébrique complexe
(et ne permet donc pas de prouver le théorème de Weil).

\cite{Serre-1960} a ensuite élaboré cet argument en dimension
supérieure pour démontrer un \emph{analogue kählérien.}
Alors, $V$ est une variété projective compacte
munie d'un endomorphisme $f\colon V\to V$
tel que l'image inverse $f^*H$ d'une section hyperplane~$H$
soit numériquement équivalente à~$qH$, pour un nombre réel $q>0$.
Serre démontre alors que les valeurs propres de $f^*$
agissant sur l'espace de cohomologie $H^i(V)$
sont de module $q^{i/2}$. Sa preuve utilise de façon cruciale
le théorème de Lefschetz difficile et la signature
de la forme d'intersection en restriction 
à la cohomologie primitive. Pour la transposer
aux variétés sur les corps finis, 
il faut disposer de l'analogue de ces deux résultats
pour une cohomologie de Weil. 
Si le théorème de Lefschetz difficile est actuellement
connu (mais comme \emph{conséquence}  des conjectures de Weil),
la signature des formes d'intersection est une question largement ouverte.

Dans une lettre à Serre datée du 27 août 1965, 
Grothendieck présente un faisceau de conjectures,
inspiré par la démonstration de Serre,
qui entraîneraient les conjectures de Weil.
Un aspect important est que ces conjectures ne font plus référence
à une cohomologie de Weil mais s'expriment uniquement
en termes des cycles algébriques sur la variété~$V$
et leurs intersections: ce sont les \emph{conjectures standard},
et je renvoie à \cite{Grothendieck-1969} et \cite{Kleiman-1968,Kleiman-1994}
pour une présentation.

\subsection{Le cas général}

La démonstration générale de l'hypothèse de Riemann pour les variétés
sur un corps fini est due à \cite{Deligne-1974a}.
Il m'est impossible de décrire cette preuve
qui repose
sur la machinerie de la cohomologie étale $\ell$-adique,
mais d'une façon moins « statique » que les approches antérieures.
En particulier, intervient l'idée d'exploiter 
qu'une famille de sections hyperplanes donne lieu à une 
grosse monodromie.

Deligne généralisera ses arguments pour établir
une propriété de stabilité fondamentale des \emph{poids}
en cohomologie étale, voir \cite{Deligne-1980}.

\section{Applications}

Il est probablement impossible de faire la liste de toutes les
applications des conjectures de Weil. Je vais donner quelques
indications dans trois directions différentes : arithmétique,
géométrie et théorie des modèles.

\subsection{Arithmétique}

Comme je l'ai déjà évoqué, l'hypothèse de Riemann
sur les corps finis est intimement liée à de bonnes
majorations de sommes d'exponentielles.
Dans le cas des courbes, l'article de \citet{Weil-1948a}
détaille ce que l'on peut obtenir.
En dimension supérieure, \citet[\S8.4]{Deligne-1974a} expose une
majoration générale, due à Bombieri dont 
voici une simplification lorsque $k=\F_p$: si $f\in \Z[T_1,\dots,T_n]$
est un polynôme de degré~$d$ dont la partie homogène modulo~$p$ définit
une hypersurface lisse de degré~$d$ de $\P_{n-1}$, alors
\[ \abs{\sum_{x\in(\Z/p\Z)^n} \exp\big( 2i\pi \frac{f(x)}p \big)}
\leq (d-1)^n p^{n/2}. \]
Le livre de \cite{Katz-1980} présente une introduction à l'utilisation
de la cohomologie étale pour les sommes d'exponentielles.

Une autre application importante est liée à la théorie des formes modulaires.
Le produit infini
\[ \Delta(z) = q \prod_{n=1}^\infty (1-q^n)^{24}, \qquad q = \exp(2i\pi z) \]
définit une forme modulaire de poids~$12$ pour le groupe $\SL(2,\Z)$;
explicitement, on a 
\[ \Delta\big( \frac{a z+b}{cz+d}\big) = (cz+d)^{12} \Delta(z).  \]
\citet{Ramanujan-1916} avait initié l'étude 
des coefficients $\tau(n)$ du développement en série
\[ \Delta(z) = \sum_{n=1}^\infty \tau(n) q^n \]
de~$\Delta$ en conjecturant leur multiplicativité
($\tau(m)\tau(n)=\tau(mn)$ si $m$ et $n$ sont premiers entre eux)
et, plus généralement, la décomposition en « produit eulérien »
\[ \sum_{n=1}^\infty \tau(n) n^{-s} = \prod_{\text{$p$ premier}} \dfrac1{1-\tau(p)p^{-s}+p^{11-2s}}, \]
ainsi que l'inégalité 
\[ \abs{\tau(p)} < 2 p^{11/2}, \]
pour tout nombre premier~$p$.
La décomposition en produit eulérien est un théorème de \cite{Mordell-1917},
mais cette dernière majoration devra attendre que \cite{Deligne-1971},
généralisant des idées de Kuga et Shimura, 
la déduise de l'hypothèse de Riemann. 
La raison d'être de cette inégalité est que dans la factorisation 
\[ 1 - \tau(p)p^{-s}+p^{11-2s} = (1 - \alpha_p p^{-s})  (1-\beta_p p^{-s}), \]
les racines $\alpha_p,\beta_p$ du polynôme $1-\tau(p) x+p^{11}x^2$, 
sont des valeurs propres de l'endomorphisme de Frobenius $\phi_V$
agissant sur un espace d'homologie $H_11(V)$ d'une variété lisse. 
Les détails sont évidemment un peu plus compliqués…

Citons une application de ces majorations de Ramanujan
par \cite{LubotzkyPhillipsSarnak-1988} à une construction
optimale de graphes expanseurs.

\subsection{Géométrie}

La première application géométrique de la démonstration
de l'hypothèse de Riemann est peut-être la démonstration
de l'analogue en cohomologie étale
du théorème de Lefschetz difficile
\citep[Théorème 4.1.1]{Deligne-1980}.
C'est un énoncé pour les variétés algébriques sur un corps 
algébriquement clos, si bien qu'il convient peut-être d'expliquer
comment il peut bien être lié à un énoncé sur le nombre
de points de variétés sur un corps fini.

De fait, il y a plusieurs procédés de réduction pour étudier
la géométrie des variétés sur un corps~$k$.
On peut notamment observer que la définition de ces variétés 
et des divers objets géométriques qui interviennent dans un énoncé
donné ne font intervenir qu'un nombre fini de polynômes;
voyant  les coefficients de ces polynômes comme des
indéterminées, la situation apparaît comme une spécialisation
d'un énoncé « en famille ». Le principe est alors
que souvent, la validité du théorème en vue est uniforme dans la famille,
si bien qu'il suffirait de le prouver pour de bonnes spécialisations.

Une extension de ce principe est au cœur de la démonstration
du « théorème de décomposition » de \cite{BeilinsonBernsteinDeligneGabber-2018}.
Une application récente de ce principe est due à \citet{Batyrev-1999}:
si deux variétés algébriques complexes
de type Calabi-Yau (projectives, lisses, à fibré canonique trivial)
sont birationnelles, elles ont mêmes nombres de Betti.
Citons aussi l'article d'exposition de \cite{Serre-2009a}.

Mentionnons enfin
que \citet{KatzMessing-1974} ont déduit du théorème de Deligne
que le théorème de Lefschetz difficile valait pour \emph{toute}
cohomologie de Weil. Ils ont aussi prouvé
que le polynôme caractéristique des endomorphismes de Frobenius de variétés
projectives lisses ne dépendait pas du choix de cette cohomologie de Weil.
Ils obtiennent en particulier que la dimension de ces espaces
(dans le cas projectif lisse)
ne dépend pas du choix de la cohomologie de Weil.

\subsection{Théorie des modèles}

\citet{Ax-1968} a déduit des estimées 
de \cite{LangWeil-1954} en direction l'hypothèse de Riemann
une caractérisation remarquable des corps finis de « grand cardinal ».

Il se place en théorie des modèles c'est-à-dire qu'il s'intéresse
aux \emph{formules} qui sont vraies dans un corps fini.
Les formules qu'il considère utilisent le « langage des anneaux »:
elles sont écrites à l'aide des symboles $+$ (addition), $-$ (soustraction), $\cdot$ (multiplication),
$0$ et $1$ (zéro et un), 
de symboles de variables et de quantificateurs $\exists$ 
(il existe), $\forall$ (pour tout), $\neg$ (négation). 
Par exemple, $\exists x,x^2+1=0$ est une telle formule,
à condition d'écrire $x^2=x\cdot x$, 
mais $\exists x\exists n,x^n=2$ n'en est pas une car on ne dispose
pas de symbole de puissance entière.
Si une telle formule $\phi$ n'a pas de variable libre et si 
$k$ est un corps, on dit qu'elle est vérifiée dans ce corps
(on note $k\models \phi$)  si son interprétation évidente
est vraie 
lorsqu'on lit les quantificateurs $\exists x$ comme « il existe $x\in k$ »,
etc. 
Par exemple, si $p$ est un nombre premier,
la formule « $\exists x, x^2+1=0$ » est vérifiée dans $\Z/p\Z$
si et seulement si $p$ n'est pas congru à~$3$ modulo~$4$;
elle est vérifiée dans $\C$, mais pas dans~$\R$ ni dans $\Q(i\sqrt 3)$.

Avec un peu de  travail en géométrie algébrique,
on peut écrire à l'aide d'une telle
formule le contexte de l'estimation de \cite{LangWeil-1954}:
il s'agit d'exprimer l'expression « soit $V$ une variété algébrique
affine géométriquement intègre et de dimension~$d$ » en termes concrets:
$V$ est définie par une famille finie $(f_1,\dots,f_r)$ de polynômes
de degrés $\leq D$, 
et il s'agit d'exprimer 
que $V$ est géométriquement intègre et de dimension~$d$
par une formule $\phi_{d,r,D}$
du langage des anneaux en les coefficients de ces polynômes.
Les entiers $d, r, D$ étant fixés, le théorème de \cite{LangWeil-1954} 
fournit alors une estimation uniforme 
\[ \abs{\Card(V(k)) - q^d} \leq c_{d,r,D} q^{d-\frac12} \]
dès que $k$ est un corps fini de cardinal~$q$.
En particulier, si $q> c_{d,r,D}^2$, on en déduit
que $V(k)\neq\emptyset$.
En particulier, la formule $\phi'_{d,r,D}$
qui exprime que pour toute telle~$V$, l'ensemble $V(k)$ n'est pas vide
est vraie pour tout corps  fini~$k$ de cardinal assez grand.

Le théorème principal d'\cite{Ax-1968}
est qu'à l'inverse, une formule $\phi$ est vraie
dans tout corps fini de cardinal assez grand si et seulement
si elle se déduit des formules $\phi'_{d,r,D}$.

\section{Au delà des bornes de Weil}

Malgré leur importance, il n'est peut-être pas inutile de signaler
que les majorations qu'entraînent les conjectures de Weil
pour le nombre de solutions d'équations dans les corps finis
sont loin d'épuiser ce qu'on voudrait savoir de ce nombre.

Je donne deux exemples.

a) Tout d'abord, déjà dans le cas d'une courbe projective lisse de genre~$g$,
la majoration $\abs{N(V)-q-1}\leq 2g\sqrt q$
n'est pas forcément applicable à une question donnée,
notamment lorsqu'il s'agit de faire varier~$V$.

Si l'on cherche une minoration de $N(V)$, l'inégalité 
\[ N(V) \geq q+1 - 2g \sqrt q \]
est par exemple inutile lorsque $g$ est trop grand.
De fait, il est possible d'avoir $N(V)=0$.

Étant donnée une suite de courbes $(V_g)$ de genres tendant vers l'infini,
définies sur un même corps fini~$k$ à $q$ éléments,
la majoration de Weil entraîne
\[ \varlimsup \frac1{2g} \Card(V_g(k)) \leq \sqrt q. \]
Lorsque $q$ est un carré parfait, on sait produire
une suite de courbes qui réalise cette limite supérieure;
dans le cas contraire, on ne sait pas s'il en existe.

b) Plutôt que des majorations/minorations, on peut avoir besoin
de congruences pour l'entier $N(V)$,
dans l'esprit du classique théorème de Chevalley--Warning,
voir \cite{Chevalley-1935b,Warning-1935}.

La cohomologie étale fournit en principe
de telles congruences modulo les nombres premiers~$\ell$  
qui ne divisent pas la caractéristique~$p$ de~$k$,
mais elle en donne peu car ces congruences sont rares.
Dans le cas des formes modulaires, 
voir par exemple \cite{Serre-1973}.

On dispose en revanche de congruences $p$-adiques,
démontrées soit par des méthodes de type Dwork, 
par exemple \citep{Ax-1964},
soit par l'étude précise de la cohomologie cristalline,
cf. \citep{Mazur-1972,Mazur-1973,Mazur-1975}.
Un exemple récent d'application de ces idées à l'\emph{existence}
de points rationnels est le théorème d'\cite{Esnault-2003}:
si $V$ est une variété projective lisse, géométriquement
rationnellement connexe, sur un corps fini~$k$ de cardinal~$q$, 
alors $V(k)$ n'est pas vide; plus précisément,
$\Card(V(k))\equiv 1 \pmod q$. 
(Cela s'applique en particulier aux variétés de Fano.)

\bibliographystyle{mynat}
\itemsep 0pt plus 1pt
\nocite{OortSchappacher-2016,Milne-2016}
\bibliography{aclab,Weil-Zotero}

\providecommand{\noopsort}[1]{}\providecommand{\url}[1]{\textit{#1}}
\begin{thebibliography}{56}
\ProvideTextCommand{\guillemotleft}{OT1}{%
  \leavevmode\raise .27ex\hbox{$\scriptscriptstyle\ll$}}
\ProvideTextCommand{\guillemotright}{OT1}{%
  \leavevmode\raise .27ex\hbox{$\scriptscriptstyle\gg$}}
\newcommand{\enquote}[1]{\og #1\fg}
\providecommand{\og}{``}\providecommand{\fg}{''}
\expandafter\ifx\csname natexlab\endcsname\relax\def\natexlab#1{#1}\fi
\expandafter\ifx\csname url\endcsname\relax
  \def\url#1{\texttt{#1}}\fi
\expandafter\ifx\csname urlprefix\endcsname\relax\def\urlprefix{\textsc{url:}
  }\fi
\newcommand{\Capitalize}[1]{\uppercase{#1}}
\newcommand{\capitalize}[1]{\expandafter\Capitalize#1}
\providecommand{\eprint}[2][]{\url{#2}}

\bibitem[{Artin(1924{\natexlab{\emph{a}}})}]{Artin-1924a}
E.~\textsc{Artin} (1924{\natexlab{\emph{a}}}), \enquote{{Quadratische K\"orper
  im Gebiete der h\"oheren Kongruenzen. I.}} \emph{Mathematische Zeitschrift},
  \textbf{19}~(1), \bblpp{} 153--206.

\bibitem[{Artin(1924{\natexlab{\emph{b}}})}]{Artin-1924}
E.~\textsc{Artin} (1924{\natexlab{\emph{b}}}), \enquote{{Quadratische K\"orper
  im Gebiete der h\"oheren Kongruenzen. II.}} \emph{Mathematische Zeitschrift},
  \textbf{19}~(1), \bblpp{} 207--246.

\bibitem[{Ax(1964)}]{Ax-1964}
J.~\textsc{Ax} (1964), \enquote{Zeroes of polynomials over finite fields}.
  \emph{American Journal of Mathematics}, \textbf{86}, \bblpp{} 255--261.

\bibitem[{Ax(1968)}]{Ax-1968}
J.~\textsc{Ax} (1968), \enquote{The {{Elementary Theory}} of {{Finite
  Fields}}}. \emph{The Annals of Mathematics}, \textbf{88}~(2), \bblp{} 239.

\bibitem[{Batyrev(1999)}]{Batyrev-1999}
V.~V. \textsc{Batyrev} (1999), \enquote{Stringy {{Hodge}} numbers of varieties
  with {{Gorenstein}} canonical singularities}. \emph{Integrable {{Systems}}
  and {{Algebraic Geometry}}}, \bblpp{} 1--32, {World Sci. Publ.},
  {Kobe/Kyoto}.

\bibitem[{Be{\u \i}linson \emph{\bbletal{}}(2018 [1982])Be{\u \i}linson,
  Bernstein, Deligne \& Gabber}]{BeilinsonBernsteinDeligneGabber-2018}
A.~A. \textsc{Be{\u \i}linson}, J.~\textsc{Bernstein}, P.~\textsc{Deligne} \&
  O.~\textsc{Gabber} (2018 [1982]), \enquote{Faisceaux pervers}. \emph{Analyse
  et Topologie Sur Les Espaces Singuliers, {{I}} ({{Luminy}}, 1981)},
  Ast\'erisque~\textbf{100}, \bblpp{} 5--171, {Soc. Math. France}, {Paris},
  2\textsuperscript{e} \bbledition{}.

\bibitem[{Bombieri(1974)}]{Bombieri-1974}
E.~\textsc{Bombieri} (1974), \enquote{Counting points on curves over finite
  fields}. \emph{S\'eminaire {{Bourbaki}} Vol. 1972/73 {{Expos\'es}}
  418\textendash 435}, Lecture {{Notes}} in {{Mathematics}}, \bblpp{} 234--241,
  {Springer}, {Berlin, Heidelberg}.

\bibitem[{Borel(1894)}]{Borel-1894}
{\'E}.~\textsc{Borel} (1894), \enquote{Sur une application d'un th\'eor\`eme de
  {{M}}. {{Hadamard}}}. \emph{Bull. Sc. Math.}, \textbf{18}, \bblpp{} 22--25.

\bibitem[{Chevalley(1935)}]{Chevalley-1935b}
C.~\textsc{Chevalley} (1935), \enquote{{D\'emonstration d'une hypoth\`ese de M.
  Artin}}. \emph{Abhandlungen aus dem Mathematischen Seminar der Universit\"at
  Hamburg}, \textbf{11}~(1), \bblpp{} 73--75.

\bibitem[{Chowla(1949)}]{Chowla-1949}
S.~\textsc{Chowla} (1949), \enquote{The {{Last Entry}} in {{Gauss}}'s
  {{Diary}}}. \emph{Proceedings of the National Academy of Sciences of the
  United States of America}, \textbf{35}~(5), \bblpp{} 244--246.

\bibitem[{Davenport(1933)}]{Davenport-1933}
H.~\textsc{Davenport} (1933), \enquote{{On certain exponential sums}}.
  \emph{Journal f\"ur die reine und angewandte Mathematik}, \textbf{169},
  \bblpp{} 158--176.

\bibitem[{Dedekind(1857)}]{Dedekind-1857}
R.~\textsc{Dedekind} (1857), \enquote{{Abriss einer Theorie der h\"ohern
  Congruenzen in Bezug auf einen reellen Primzahl-Modulus.}} \emph{Journal
  f\"ur Mathematik}, \textbf{54}~(1), \bblpp{} 1--26.

\bibitem[{Deligne(1971)}]{Deligne-1971}
P.~\textsc{Deligne} (1971), \enquote{Formes modulaires et repr\'esentations
  $\ell$-adiques}. \emph{S\'eminaire {{Bourbaki}} 1968/69}, Lecture {{Notes}}
  in {{Math}}.~\textbf{175}, \bblpp{} Expos\'e 355, 139--172, {Springer}.

\bibitem[{Deligne(1974)}]{Deligne-1974a}
P.~\textsc{Deligne} (1974), \enquote{{La conjecture de Weil. I}}.
  \emph{Publications math\'ematiques de l'IH\'ES}, \textbf{43}~(1), \bblpp{}
  273--307.

\bibitem[{Deligne(1980)}]{Deligne-1980}
P.~\textsc{Deligne} (1980), \enquote{{La Conjecture de Weil. II}}.
  \emph{Publications math\'ematiques de l'IH\'ES}, \textbf{52}~(1), \bblpp{}
  137--252.

\bibitem[{Dwork(1960)}]{Dwork-1960a}
B.~\textsc{Dwork} (1960), \enquote{On the rationality of the zeta function of
  an algebraic variety}. \emph{American Journal of Mathematics},
  \textbf{82}~(3), \bblpp{} 631--648.

\bibitem[{Esnault(2003)}]{Esnault-2003}
H.~\textsc{Esnault} (2003), \enquote{Varieties over a finite field with trivial
  {{Chow}} group of 0-cycles have a rational point}. \emph{Inventiones
  mathematicae}, \textbf{151}~(1), \bblpp{} 187--191.

\bibitem[{Gauss(1863)}]{Gauss-1863b}
C.~F. \textsc{Gauss} (1863), \emph{{Werke II. H\"ohere Arithmetik}},
  {K\"oniglichen Gesellschaft der Wissenschaften zu G\"ottingen},
  {G\"ottingen}.

\bibitem[{Gauss(1863 [1801])}]{Gauss-1863a}
C.~F. \textsc{Gauss} (1863 [1801]), \emph{Werke {{I}}. {{Disquisitiones}}
  Arithmeticae}, {K\"oniglichen Gesellschaft der Wissenschaften zu
  G\"ottingen}, {G\"ottingen}.

\bibitem[{Grothendieck(1958)}]{Grothendieck-1958b}
A.~\textsc{Grothendieck} (1958), \enquote{{Sur une note de Mattuck-Tate.}}
  \emph{Journal f\"ur die reine und angewandte Mathematik}, \textbf{200},
  \bblpp{} 208--215.

\bibitem[{Grothendieck(1969)}]{Grothendieck-1969}
A.~\textsc{Grothendieck} (1969), \enquote{Standard conjectures on algebraic
  cycles}. \emph{Algebraic {{Geometry}}. 1968, 193-199 (1969).}, \bblpp{}
  193--199, {Bombay, 1968}.

\bibitem[{Hasse(1936)}]{Hasse-1936}
H.~\textsc{Hasse} (1936), \enquote{Zur {{Theorie}} der abstrakten elliptischen
  {{Funktionenk\"orper III}}. {{Die Struktur}} des {{Meromorphismenrings}}.
  {{Die Riemannsche Vermutung}}.} \emph{jreia}, \textbf{175}~(4), \bblpp{}
  193--208.

\bibitem[{Hasse \& Davenport(1935)}]{HasseDavenport-1935}
H.~\textsc{Hasse} \& H.~\textsc{Davenport} (1935), \enquote{{Die Nullstellen
  der Kongruenzzetafunktionen in gewissen zyklischen F\"allen.}} \emph{Journal
  f\"ur die reine und angewandte Mathematik}, \textbf{172}, \bblpp{} 151--182.

\bibitem[{Herglotz(1921)}]{Herglotz-1921}
G.~\textsc{Herglotz} (1921), \enquote{{Zur letzten Eintragung im Gau\ss schen
  Tagebuch}}. \emph{Ber. Math. Phys. Kl. S\"achs. Akad. Wiss. Leipzig},
  \textbf{73}, \bblpp{} 271--276.

\bibitem[{Katz(1980)}]{Katz-1980}
N.~M. \textsc{Katz} (1980), \emph{{Sommes exponentielles. Cours \`a Orsay,
  automne 1979. Redige par Gerard Laumon, preface par Luc Illusie}},
  {Ast\'erisque}~\textbf{79}, {Soci\'et\'e Math\'ematique de France (SMF),
  Paris}.

\bibitem[{Katz \& Messing(1974)}]{KatzMessing-1974}
N.~M. \textsc{Katz} \& W.~\textsc{Messing} (1974), \enquote{Some consequences
  of the {{Riemann}} hypothesis for varieties over finite fields}.
  \emph{Inventiones mathematicae}, \textbf{23}~(1), \bblpp{} 73--77.

\bibitem[{Kleiman(1968)}]{Kleiman-1968}
S.~L. \textsc{Kleiman} (1968), \enquote{Algebraic cycles and the {{Weil}}
  conjectures}. \emph{Dix Expos\'es Sur La Cohomologie Des Sch\'emas}, Adv.
  {{Stud}}. {{Pure Math}}.~\textbf{3}, \bblpp{} 359--386, {North-Holland,
  Amsterdam}.

\bibitem[{Kleiman(1994)}]{Kleiman-1994}
S.~L. \textsc{Kleiman} (1994), \enquote{The standard conjectures}.
  \emph{Motives. {{Proceedings}} of the Summer Research Conference on Motives,
  Held at the {{University}} of {{Washington}}, {{Seattle}}, {{WA}}, {{USA}},
  {{July}} 20-{{August}} 2, 1991}, \bblpp{} 3--20, {American Mathematical
  Society}, {Providence, RI}.

\bibitem[{Klein(1903)}]{Klein-1903}
F.~\textsc{Klein} (1903), \enquote{{Gau\ss ' wissenschaftliches Tagebuch
  1796\textendash 1814}}. \emph{Mathematische Annalen}, \textbf{57}~(1),
  \bblpp{} 1--34.

\bibitem[{Lang \& Weil(1954)}]{LangWeil-1954}
S.~\textsc{Lang} \& A.~\textsc{Weil} (1954), \enquote{Number of points of
  varieties in finite fields}. \emph{American Journal of Mathematics},
  \textbf{76}, \bblpp{} 819--827.

\bibitem[{Lubotzky \emph{\bbletal{}}(1988)Lubotzky, Phillips \&
  Sarnak}]{LubotzkyPhillipsSarnak-1988}
A.~\textsc{Lubotzky}, R.~\textsc{Phillips} \& P.~\textsc{Sarnak} (1988),
  \enquote{Ramanujan graphs}. \emph{Combinatorica}, \textbf{8}~(3), \bblpp{}
  261--277.

\bibitem[{Mattuck \& Tate(1958)}]{MattuckTate-1958}
A.~\textsc{Mattuck} \& J.~\textsc{Tate} (1958), \enquote{On the inequality of
  {{Castelnuovo-Severi}}: {{To Emil Artin}} on his 60th birthday}.
  \emph{Abhandlungen aus dem Mathematischen Seminar der Universit\"at Hamburg},
  \textbf{22}~(1), \bblpp{} 295--299.

\bibitem[{Mazur(1972)}]{Mazur-1972}
B.~\textsc{Mazur} (1972), \enquote{Frobenius and the {{Hodge}} filtration}.
  \emph{Bulletin of the American Mathematical Society}, \textbf{78}~(5),
  \bblpp{} 653--667.

\bibitem[{Mazur(1973)}]{Mazur-1973}
B.~\textsc{Mazur} (1973), \enquote{Frobenius and the {{Hodge}} filtration
  (estimates)}. \emph{Annals of Mathematics. Second Series}, \textbf{98},
  \bblpp{} 58--95.

\bibitem[{Mazur(1975)}]{Mazur-1975}
B.~\textsc{Mazur} (1975), \enquote{Eigenvalues of {{Frobenius}} acting on
  algebraic varieties over finite fields}. \emph{Proc. {{Symp}}. {{Pure
  Math}}.}, ~\textbf{29}, \bblpp{} 231--261, {American Mathematical Society},
  {Arcata}.

\bibitem[{Milne(2016)}]{Milne-2016}
J.~S. \textsc{Milne} (2016), \enquote{The {{Riemann}} hypothesis over finite
  fields: From {{Weil}} to the present day}. \emph{Notices of the International
  Congress of Chinese Mathematicians}, \textbf{4}~(2), \bblpp{} 14--52.

\bibitem[{Mordell(1917)}]{Mordell-1917}
L.~J. \textsc{Mordell} (1917), \enquote{On {{Mr}}. {{Ramanujan}}'s empirical
  expansions of modular functions}. \emph{Proceedings of the Cambridge
  Philosophical Society}, \textbf{19}, \bblpp{} 117--124.

\bibitem[{Oort \& Schappacher(2016)}]{OortSchappacher-2016}
F.~\textsc{Oort} \& N.~\textsc{Schappacher} (2016), \enquote{Early {{History}}
  of the {{Riemann Hypothesis}} in {{Positive Characteristic}}}. \emph{The
  {{Legacy}} of {{Bernhard Riemann After One Hundred}} and {{Fifty Years}}},
  \bbledby{} L.~\textsc{Ji}, F.~\textsc{Oort} \& S.-T. \textsc{Yau},
  \textbf{35}, \bblpp{} 595--631, {Higher Education Press and International
  Press}, {Beijing\textendash Boston}.

\bibitem[{Ramanujan(1916)}]{Ramanujan-1916}
S.~\textsc{Ramanujan} (1916), \enquote{On certain arithmetical functions}.
  \emph{Transactions of the Cambridge Philosophical Society}, \textbf{22},
  \bblpp{} 159--184.

\bibitem[{Roquette(2018)}]{Roquette-2018}
P.~\textsc{Roquette} (2018), \emph{The {{Riemann Hypothesis}} in
  {{Characteristic}} p in {{Historical Perspective}}}, Lecture {{Notes}} in
  {{Mathematics}}~\textbf{2222}, {Springer International Publishing}, {Cham}.

\bibitem[{Schmidt(1931)}]{Schmidt-1931}
F.~K. \textsc{Schmidt} (1931), \enquote{{Analytische Zahlentheorie in K\"orpern
  der Charakteristikp}}. \emph{Mathematische Zeitschrift}, \textbf{33}~(1),
  \bblpp{} 1--32.

\bibitem[{Schmidt(1976)}]{Schmidt-1976}
W.~M. \textsc{Schmidt} (1976), \emph{Equations over {{Finite Fields An
  Elementary Approach}}}, Lecture {{Notes}} in {{Mathematics}}~\textbf{536},
  {Springer Berlin Heidelberg}, {Berlin, Heidelberg}.

\bibitem[{Serre(1960)}]{Serre-1960}
J.-P. \textsc{Serre} (1960), \enquote{Analogues k\"ahl\'eriens de certaines
  conjectures de {{Weil}}}. \emph{Annals of Mathematics. Second Series},
  \textbf{71}, \bblpp{} 392--394.

\bibitem[{Serre(1973)}]{Serre-1973}
J.-P. \textsc{Serre} (1973), \enquote{Congruences et formes modulaires
  (d'apr\`es {{H}}.{{P}}.{{F}}. {{Swinnerton-Dyer}})}. \emph{S\'eminaire
  {{Bourbaki}} 1971/72}, Lecture {{Notes}} in {{Mathematics}}~\textbf{317},
  \bblpp{} Expos\'e 416, 319--338.

\bibitem[{Serre(2009)}]{Serre-2009a}
J.-P. \textsc{Serre} (2009), \enquote{How to use finite fields for problems
  concerning infinite fields}. \emph{Arithmetic, Geometry, Cryptography and
  Coding Theory}, Contemp. {{Math}}.~\textbf{487}, \bblpp{} 183--193, {Amer.
  Math. Soc., Providence, RI}.

\bibitem[{Stepanov(1969)}]{Stepanov-1969}
S.~A. \textsc{Stepanov} (1969), \enquote{{On the number of points of a
  hyperelliptic curve over a finite prime field}}. \emph{Izvestiya Akademii
  Nauk SSSR. Seriya Matematicheskaya}, \textbf{33}, \bblpp{} 1171--1181.

\bibitem[{Ullrich(2000)}]{Ullrich-2000a}
P.~\textsc{Ullrich} (2000), \enquote{{Emil Artin's unpublished generalization
  of his dissertation}}. \emph{Mitteilungen der Mathematischen Gesellschaft in
  Hamburg}, \textbf{19}, \bblpp{} 173--194.

\bibitem[{Warning(1935)}]{Warning-1935}
E.~\textsc{Warning} (1935), \enquote{{Bemerkung zur vorstehenden Arbeit von
  Herrn Chevalley}}. \emph{Abhandlungen aus dem Mathematischen Seminar der
  Universit\"at Hamburg}, \textbf{11}~(1), \bblpp{} 76--83.

\bibitem[{Weil(1940)}]{Weil-1940}
A.~\textsc{Weil} (1940), \enquote{{Sur les fonctions alg\'ebriques \`a corps de
  constantes fini}}. \emph{Comptes Rendus Hebdomadaires des S\'eances de
  l'Acad\'emie des Sciences, Paris}, \textbf{210}, \bblpp{} 592--594.

\bibitem[{Weil(1941)}]{Weil-1941}
A.~\textsc{Weil} (1941), \enquote{On the {{Riemann Hypothesis}} in
  {{Function-Fields}}}. \emph{Proceedings of the National Academy of Sciences},
  \textbf{27}~(7), \bblpp{} 345--347.

\bibitem[{Weil(1946)}]{Weil-1946}
A.~\textsc{Weil} (1946), \emph{Foundations of Algebraic Geometry}, Colloq.
  {{Publ}}., {{Am}}. {{Math}}. {{Soc}}.~\textbf{29}, {American Mathematical
  Society (AMS), Providence, RI}.

\bibitem[{Weil(1948{\natexlab{\emph{a}}})}]{Weil-1948a}
A.~\textsc{Weil} (1948{\natexlab{\emph{a}}}), \enquote{On some exponential
  sums}. \emph{Proceedings of the National Academy of Sciences of the United
  States of America}, \textbf{34}, \bblpp{} 204--207.

\bibitem[{Weil(1948{\natexlab{\emph{b}}})}]{Weil-1948b}
A.~\textsc{Weil} (1948{\natexlab{\emph{b}}}), \emph{{Sur les courbes
  alg\'ebriques et les vari\'et\'es qui s'en d\'eduisent}}, {Actualit\'es Sci.
  Ind.}~\textbf{1041}, {Hermann \& Cie.}, {Paris}.

\bibitem[{Weil(1948{\natexlab{\emph{c}}})}]{Weil-1948}
A.~\textsc{Weil} (1948{\natexlab{\emph{c}}}), \emph{Vari\'et\'es Ab\'eliennes
  et Courbes Alg\'ebriques}, Actualit\'es {{Sci}}. {{Ind}}.~\textbf{1064},
  {Hermann \& Cie., Paris}.

\bibitem[{Weil(1949)}]{Weil-1949}
A.~\textsc{Weil} (1949), \enquote{Number of solutions of equations in finite
  fields}. \emph{Bull. Amer. Math. Soc.}, \textbf{55}, \bblpp{} 397--508.

\bibitem[{Weil(1956)}]{Weil-1956a}
A.~\textsc{Weil} (1956), \enquote{Abstract versus classical algebraic
  geometry}. \emph{Proceedings of the {{International Congress}} of
  {{Mathematicians}} 1954. {{Amsterdam}}, {{September}} 2\textendash 9.
  {{Vol}}. {{III}}. {{Stated}} Addresses in Sections. {{Symposia}}}, {Erven P.
  Noordhoff N. V.; Amsterdam: North-Holland Publishing Co.}, {Groningen}.

\end{thebibliography}

\end{document}